\documentclass[12pt]{amsart}
\usepackage{verbatim}
\usepackage{amsmath}
\usepackage{amsthm}
\usepackage{amssymb}
\usepackage{enumerate}

\newif\ifdeveloping

\newtheorem{theorem}{Theorem}

\newtheorem{problem}[theorem]{Problem}

\theoremstyle{definition}

\newtheorem{definition}[theorem]{Definition}

\theoremstyle{remark}
  
{}

\ifdeveloping
\usepackage[notref,notcite]{showkeys}
\fi

\newcommand{\prtime}{{\count0=\time\divide\count0 by 60
\count1=-\count0\multiply\count1 by 60 \advance\count1 by \time
\the\count0:\the\count1} }

\def\myheads#1;#2;{
\pagestyle{myheadings} \markboth{{\sc\hfill
#1\hfill\protect\makebox[0cm][r]{\rm\today; \prtime}}}
{{\sc\protect\makebox[0cm][l]{\rm\today;\ \prtime}\hfill #2\hfill}}
\thispagestyle{myheadings} }

\newcommand{\acal}{{\mathcal A}}

\newcommand{\fcal}{{\mathcal F}}

\newcommand{\subs}{\subset}

\newcommand{\dom}{\operatorname{dom}}

\def\<{\left\langle}
\def\>{\right\rangle}
\def\cf{\operatorname{cf}}

\def\br#1;#2;{\bigl[ {#1} \bigr]^ {#2} }
\def\bc#1;#2;{\bigl( {#1} \bigr)^ {#2} }

\theoremstyle{plain}

\begin{document}
\author[I. Juh\'asz]{Istv\'an Juh\'asz}
\address{Alfr{\'e}d R{\'e}nyi Institute of Mathematics}
\email{juhasz@renyi.hu}
\author[S. Shelah]{Saharon Shelah}
\address{Einstein Institute of Mathematics,
The Hebrew University of Jerusalem}
\email{shelah@math.huji.ac.il}
\thanks{The first author was supported
by OTKA grant no.\ 61600.}

\title[Lindel\" of]{Lindel\" of spaces of singular density}

\begin{abstract}
A cardinal $\lambda$ is called $\omega$-inaccessible if for all $\mu
< \lambda$ we have $\mu^\omega < \lambda.$ We show that for every
$\omega$-inaccessible cardinal $\lambda$ there is a CCC (hence
cardinality and cofinality preserving) forcing that adds a
hereditarily Lindel\" of regular space of density $\lambda.$ This
extends an analogous earlier result of ours that only worked for
regular $\lambda$.
\end{abstract}



\maketitle

In \cite{JSh} we have shown that for any cardinal $\lambda$ a
natural CCC forcing notion adds a hereditarily Lindel\"of
0-dimensional Hausdorff topology on $\lambda$ that makes the
resulting space $X_\lambda\,$ left-separated in its natural
well-ordering. It was also shown there that the density
$d(X_\lambda) = \cf(\lambda)$, hence if $\lambda$ is regular then
$d(X_\lambda) = \lambda$. The aim of this paper is to show that a
suitable extension of the construction given in \cite{JSh} enables
us to generalize this to many singular cardinals as well.

Note that the existence of an L-space, that we now know is provable
in ZFC (see \cite{M}), is equivalent to the existence of a
hereditarily Lindel\" of regular space of density $\omega_1.$ Since
the cardinality of a hereditarily Lindel\" of $T_2$ space is at most
continuum, just in ZFC we cannot replace in this $\omega_1$ with
anything bigger. The following problem however, that is left open by
our subsequent result, can be raised naturally.

\begin{problem}
Assume that $\omega_1 < \lambda \le \mathfrak{c}.$ Does there exist
then a hereditarily Lindel\" of regular space of density
$\lambda\,$?
\end{problem}

We should emphasize that this problem is open for all cardinals
$\lambda$, regular or singular, in particular for  $\lambda =
\omega_2.$

Before describing our new construction, let us recall that the one
given in \cite{JSh} is based on simultaneously and generically
``splitting into two" the complements $\lambda \setminus \alpha$ for
all proper initial segments $\alpha$ of $\lambda.$ The novelty in
the construction to be given is that we shall perform the same
simultaneous splitting for the complements of the members of a
family $\mathcal{A}$ of subsets of $\lambda$ that is, at least when
$\lambda$ is singular, much larger than the family of its proper
initial segments (that is just $\lambda$ if we are considering von
Neumann ordinals). The following definition serves to describe the
properties of such a family of subsets of $\lambda.$

\begin{definition}\label{df:A}
Let $\lambda$ be an infinite cardinal. A family $\mathcal{A}$ of
{\em proper} subsets of $\lambda$ is said to be {\em good} over
$\lambda$ if it satisfies properties (i)-(iii) below:
\begin{itemize}
\item[(i)]
$\lambda \subs \mathcal{A}$\, that is, all proper initial segments
of $\lambda$ belong to $\mathcal{A}$;
\item[(ii)]
for every subset $S \subs \lambda$ with $|S| < \lambda$ there is $A
\in \mathcal{A}$ with $S \subs A$;
\item[(iii)]
for every subset $S \subs \lambda$ with $|S| = \omega_1$ there is $T
\in [S]^{\omega_1}$ such that if $A \in \mathcal{A}$ then either $|A
\cap T| \le \omega$ or $T \subs A.$
\end{itemize}
\end{definition}

If $\lambda$ is regular then $\mathcal{A} = \lambda$, the family of
all proper initial segments of $\lambda$, is a good family over
$\lambda$. Indeed, (i) and (ii) are obviously valid and if $S \in
[\lambda]^{\omega_1}$ then any subset $T$ of $S$ of order type
$\omega_1$ satisfies (iii). If, however, $\lambda$ is singular then
this $\mathcal{A}$ definitely does not satisfy condition (ii).
Actually, we do not know if it is provable in ZFC that for any
(singular) cardinal $\lambda$ there is a good family over $\lambda$.
But we know that they do exist if $\lambda$ is
$\omega$-inaccessible, that is $\mu^\omega < \lambda$ holds whenever
$\mu < \lambda.$

\begin{theorem}\label{tm:A}
If $\lambda$ is an $\omega$-inaccessible cardinal then there exists
a good family $\acal \subs [\lambda]^{<\lambda}$ over $\lambda$.
\end{theorem}

\begin{proof}
It is well-known that there is a map $G : [\omega]^\omega
\rightarrow \omega$ with the property that for every $a \in [
\omega]^\omega$ we have $G \big[[a]^\omega \big] = \omega.$ In other
words: we may color the infinite subsets of $\omega$ with countably
many colors so that on the subsets of any infinite set all the
colors are picked up. Such a coloring may be constructed by a simple
transfinite recursion.

Next we fix a {\em maximal almost disjoint} family $\fcal$ of
subsets of order type $\omega$ of our underlying set $\lambda$ and
then we ``transfer" the coloring $G$ to each member $F$ of $\fcal.$
More precisely, this means that for every $F \in \fcal$ we fix a map
$G_F : [F]^\omega \rightarrow F$ such that $G_F \big[[a]^\omega
\big] = F$ whenever $a \in [F]^\omega.$ Then we ``fit together"
these colorings $G_F$ to obtain a coloring $H : [\lambda]^\omega
\rightarrow \lambda$ of all countable subsets of $\lambda$ as
follows: For any $S \in [\lambda]^\omega$ we set $H(S) = G_F(S)$ if
there is an $F \in \fcal$ with $S \subs F$ and $H(S) = 0$ otherwise.
The coloring $H$ is well-defined because, as $\fcal$ is almost
disjoint, for every $S \in [\lambda]^\omega$ there is at most one $F
\in \fcal$ with $S \subs F.$

Now, a set $C \subs \lambda$ is called $H$-closed if for every $S
\in [C]^\omega$ we have $H(S) \subs C.$ Clearly, for every set $A
\subs \lambda$ there is a smallest $H$-closed set including $A$ that
will be denoted by $cl_H(A)$ and is called the $H$-closure of $A.$

Let us set $A^+ = A \cup H\big[[A]^\omega\big]$ for any $A \subs
\lambda$. It is obvious that then we have $$cl_H(A) =
\bigcup_{\alpha < \omega_1} A^\alpha,$$ where the sets $A^\alpha$
are defined by the following transfinite recursion: $A^0 = A,$
$\,A^{\alpha+1} = (A^\alpha)^+,$ and $A^\alpha = \cup_{\beta <
\alpha}A^\beta$ for $\alpha$ limit. Since $$H\big[[A]^\omega\big]
\subs \bigcup \{ F \in \fcal : |F \cap A| = \omega \} \cup \{ 0
\}\,,$$ it is also obvious that we have $|A^+| \le |A|^\omega$ for
all $A \subs \lambda$ and consequently
$$|cl_H(A)| \le |A|^\omega$$ as well. In particular, $|A| < \lambda$
implies $|cl_H(A)| < \lambda$ because $\lambda$ is
$\omega$-inaccessible.

Now we claim that the family $\acal$ of all $H$-closed sets of
cardinality less than $\lambda$ is good over $\lambda$. Indeed,
first notice that because each $F \in \fcal$ has order type
$\omega$, for every set $S \in [F]^\omega$ we have $$H(S) = G_F(S) <
\sup F = \sup S,$$ implying that every initial segment $\alpha$ of
$\lambda$ is $H$-closed and so $\acal$ satisfies condition (i) of
definition \ref{df:A}. Condition (ii) is satisfied trivially.

To see (iii) we first show that there is no infinite strictly
descending sequence of $H$-closed subsets of $\lambda,$ or in other
words: the family of $H$-closed sets is well-founded with respect to
inclusion. Assume, reasoning indirectly, that $\{ C_n : n <\omega\}$
is a strictly decreasing sequence of $H$-closed sets and for each $n
< \omega$ we have $\alpha_n \in C_n \setminus C_{n+1}$. By the
maximality of $\fcal$ then there is some $F \in \fcal$ such that the
set $S = F \cap \{ \alpha_n : n < \omega \}$ is infinite. Then, for
any $k < \omega$, the set $S \cap C_k$ is also infinite and
consequently we have $$H\big[[S \cap C_k]^\omega\big] = G_F\big[[S
\cap C_k]^\omega\big] = F \subs C_k$$ because $C_k$ is $H$-closed.
But for any $k < \omega$ such that $\alpha_k \in S \subs F$ this
would imply
$$\alpha_k \in F \subs C_{k+1}\,,$$ which is clearly a contradiction.

Now let $S \subs \lambda$ with $|S| = \sigma.$ Our previous result
clearly implies that there is a set $T \in [S]^\sigma$ such that we
have $cl_H(U) = cl_H(T)$ whenever $U \subs T$ with $|U| = \sigma.$
In other words, this means that for every $H$-closed set $C$ we have
either $|C \cap T| < \sigma$ or $T \subs C.$ In particular, for
$\sigma = \omega_1$ this shows that our family $\acal$ satisfies
condition (iii) of definition \ref{df:A} as well, hence it is indeed
good over $\lambda$.
\end{proof}

Next we present our main result that, in view of theorem \ref{tm:A},
immediately implies the consistency of the existence of hereditarily
Lindel\"of regular spaces of density $\lambda$ practically for any
singular cardinal $\lambda$. (Of course, this has to be in a model
in which $\lambda \le \mathfrak{c}$.) We shall follow \cite{K} in
our notation and terminology concerning forcing.

\begin{theorem}\label{tm:d}
Let $\acal$ be a good family over $\lambda.$ Then there is a
complete (hence CCC) subforcing $\mathbb{Q}$ of the Cohen forcing
$Fn(\acal \times \lambda,\,2)$ such that in the generic extension
$V^\mathbb{Q}$ there is a hereditarily Lindel\"of 0-dimensional
Hausdorff topology $\tau$ on $\lambda$ that has density $\lambda.$
If we also have $\acal \subs [\lambda]^{< \lambda}$ (as in theorem
\ref{tm:A}) then every subset of $\lambda$ of size $< \lambda$ is
even $\tau$-nowhere dense.
\end{theorem}

\begin{proof}
We start by defining the the subforcing $\mathbb{Q}$ of $Fn(\acal
\times \lambda,\,2)$: $\mathbb{Q}$ consists of those $p \in Fn(\acal
\times \lambda,\,2)$ for which $\langle A,\alpha \rangle \in \dom p$
with $\alpha \in A $ implies $p(A,\alpha) = 0$  and $\langle A,
\gamma_A \rangle \in \dom p$ implies $p(A, \gamma_A) =1 $, where
$\gamma_A = \min(\lambda \setminus A)$. It is straight-forward to
check that $\mathbb{Q}$ is a complete suborder of $Fn(\acal \times
\lambda,\,2)$.

For any condition $p \in \mathbb{Q}$ and any set $A \in \acal$ we
define $$U^p_A = \{ \alpha : p(A,\alpha) = 1 \},$$ and if $G \subs
\mathbb{Q}$ is generic then, in $V[G]$, we set $$U_A = \bigcup \{
U^p_A : p \in G \}.$$ Next, let $U^1_A = U_A$ and $U^0_A = \lambda
\setminus A$ and $\tau$ be the topology on $\lambda$ generated by
the sets $\{ U^i_A : i < 2,\,A \in \acal \}.$ Note that then the
family $\mathcal{B} = \{ B_\varepsilon : \varepsilon \in
Fn(\acal,\,2) \}$ is a base for $\tau$, where $B_\varepsilon =
\bigcap_{A \in \dom \varepsilon}U_A^{\varepsilon(A)}$. It is clear
from the definition that each $B_\varepsilon$ is clopen, hence
$\tau$ is 0-dimensional. Now, if $\beta < \alpha < \lambda$ then we
have $\alpha \in \acal$ by (i) and hence $\beta \in \alpha \subs
U^0_\alpha$ while $\alpha = \gamma_\alpha \in U^1_\alpha$, which
shows that $\tau$ is also Hausdorff. It is also immediate from (ii)
that no set $S \in [\lambda]^{< \lambda}$ is $\tau$-dense, hence the
space $\langle \lambda, \tau \rangle$ has density $\lambda$. Indeed,
if $S \subs A \in \acal$ then we have $S \cap U^1_A = \emptyset$,
while $U^1_A \ne \emptyset.$ Thus it only remains for us to prove
that the topology $\tau$ is hereditarily Lindel\"of.

Assume, reasoning indirectly, that some condition $p \in \mathbb{Q}$
forces that $\tau$ is not hereditarily Lindel\"of, i. e. there is a
right separated $\omega_1$-sequence in $\lambda$. More precisely,
this means that there are $\mathbb{Q}$-names $\dot{s}$ and $\dot{e}$
such that $p$ forces ``$\dot{s} : \omega_1 \rightarrow \lambda,\,\,
\dot{e} : \omega_1 \rightarrow Fn(\acal,2),\,\, \dot{s}(\alpha) \in
B_{\dot{e}(\alpha)}$, and $\dot{s}(\beta) \notin
B_{\dot{e}(\alpha)}$ whenever $\alpha < \beta <\lambda$." Then, in
the ground model $V$, for each $\alpha < \omega_1$ we may pick a
condition $p_\alpha \le p$, an ordinal $\nu_\alpha < \lambda$, and a
finite function $\varepsilon_\alpha \in Fn(\acal,2)$ such that
$$p_\alpha \Vdash \dot{s}(\alpha) = \nu_\alpha \wedge \dot{e}(\alpha) = \varepsilon_\alpha.$$

Since $\mathbb{Q}$ is a complete suborder of $Fn(\acal \times
\lambda,\,2)$ it has property K, hence we may assume without any
loss of generality that the conditions $p_\alpha$ are pairwise
compatible. By extending the conditions $p_\alpha$, if necessary, we
may assume that $\dom p_\alpha = I_\alpha \times a_\alpha$ with
$I_\alpha \in [\acal]^{<\omega}$ and $a_\alpha \in
[\lambda]^{<\omega}$, moreover $\dom \varepsilon_\alpha \subs
I_\alpha$ and $\nu_\alpha \in a_\alpha$ whenever $\alpha <
\omega_1$. With an appropriate thinning out (and re-indexing) we can
achieve that if $\alpha < \beta < \omega_1$ then $$\nu_\beta \notin
a_\alpha \cup \{ \gamma_A : A \in I_\alpha \}.$$

Using standard counting and delta-system arguments, we may assume
that each $\varepsilon_\alpha$ has the same size $n < \omega$,
moreover the sets
$$\dom \varepsilon_\alpha = \{ A_{i,\alpha} : i < n\} \in [\acal]^n$$
form a delta-system, so that for some $m < n$ we have $ A_{i,\alpha}
= A_i$ if $i < m$ for all $\alpha < \omega_1$, and the families
$\{A_{m,\alpha}, ... ,A_{n-1,\alpha}\}$ are pairwise disjoint.  We
may also assume that for every $i < n$ there is a fixed value $l_i <
2$ such that $\varepsilon_\alpha(A_i) = l_i$  for all $\alpha <
\omega_1$. With a further thinning out we may achieve to have
$$\dom \varepsilon_\alpha \cap I_\beta = \{ A_i : i < m \}$$
whenever $\alpha < \beta < \omega_1$.

Finally, by property (iii) of the good family $\acal$, we may also
assume that the set $T = \{ \nu_\alpha : \alpha < \omega_1 \} \in
[\lambda]^{\omega_1} $ satisfies either $|A \cap T| \le \omega$ or
$T \subs A$ whenever $A \in \acal$.

Now, after all this thinning out, we claim that there is a countable
ordinal $\alpha > 0$ such that, for every $i < n$, if $\nu_\alpha
\in A_{i,0}$ then $l_i = 0$. Indeed, arguing indirectly, assume that
for every $0 < \alpha < \omega_1$ there is an $i_\alpha < n$ with
$\nu_\alpha \in A_{i_\alpha,0}$ and $l_{i_\alpha} = 1.$ Then there
is a fixed $j < n$ such that the set $ \{\alpha : i_\alpha = j  \}$
is uncountable and $l_j = 1$. But the first part implies $|A_{j,0}
\cap T| = \omega_1$, hence $\nu_0 \in T \subs A_{j,0} \subs
U^0_{A_{j,0}}$ that would imply $\varepsilon_0(A_{j,0}) = l_j = 0$,
a contradiction.

So, let us choose $\alpha > 0$ as in our above claim. We then define
a finite function $q \in Fn(\acal \times \lambda,\,2)$ by setting $q
\supset p_0 \cup p_\alpha$,  $$\dom q = \dom p_0 \cup \dom p_\alpha
\cup \{ \langle A_{i,0},\nu_\alpha  \rangle : m \le i < n \},$$ and
finally $$q(A_{i,0},\nu_\alpha) = l_i$$ for all $m \le i < n$. We
have $\nu_\alpha \notin a_0$, and also $A_{i,0} \notin I_\alpha$ for
$m \le i < n$ by our construction, hence this definition of $q$ is
correct. Moreover, by the above claim if $\nu_\alpha \in A_{i,0}$
then $l_i = 0$ and if $\nu_\alpha \notin A_{i,0}$ then $\nu_\alpha
\ne \gamma_{A_{i,0}}$, consequently we actually have $q \in
\mathbb{Q}$.

Let us observe, however, that we have $q(A_{i,0},\nu_\alpha) = l_i$
for all $i < n$. Indeed, if $i < m$ then this holds because
$p_\alpha(A_{i,0},\nu_\alpha) = p_\alpha(A_i,\nu_\alpha) = l_i$. But
this implies that $q \Vdash \nu_\alpha \in B_{\varepsilon_0}$ and
hence $q \Vdash \dot{s}(\alpha) \in B_{\dot{e}(0)}$ that is clearly
a contradiction because $q$ extends $p$.

Now assume that we also have $\acal \subs [\lambda]^{<\lambda}$ (in
$V$). Since $\mathbb{Q}$ is CCC, every subset of $\lambda$ in
$V^\mathbb{Q}$ is covered by a ground model set of the same size,
hence it suffices to show that any ground model member $Y$ of
$[\lambda]^{<\lambda}$ is $\tau$-nowhere dense. To see this, we
first note that it follows from a straight-forward density argument
that for every $\varepsilon \in Fn(\acal,2)$ we have
$|B_\varepsilon| = \lambda$. (Actually, this only uses the
assumption that $|\lambda \setminus \cup \acal_0| = \lambda$ for
every $\acal_0 \in [\acal]^{<\omega}$ which is weaker than $\acal
\subs [\lambda]^{<\lambda}$.)

Next, consider any set $Y \in [\lambda]^{<\lambda} \cap V$ and a
fixed $\varepsilon \in Fn(\acal,2)$. Since $\acal$ satisfies
condition (ii) of definition \ref{df:A}, we may clearly find an $A
\in \acal$ such that $Y \subs A$ and $A \notin \dom \varepsilon$.
Let $\varepsilon' = \varepsilon \cup \{ \langle A,1 \rangle \}$,
then $B_{\varepsilon'} = B_\varepsilon \cap U^1_A$ is a non-empty
open subset of $B_\varepsilon$ that is clearly disjoint from $A$ and
hence from $Y$ as well. This shows that $Y$ is indeed $\tau$-nowhere
dense.
\end{proof}

For a singular cardinal $\lambda$ of cofinality $\omega$ the results
of \cite{JSh} did imply the existence of hereditarily Lindel\" of
regular spaces of density $\lambda$, by taking the topological sum
of those of density $\lambda_n$ with $\lambda_n$ regular and
$\lambda = \sum_{n<\omega} \lambda_n$. It should be emphasized,
however, that the spaces obtained in this way clearly do not have
the stronger property we obtained in theorem \ref{tm:d} that all
subsets of size less than $\lambda$ are nowhere dense. So, we do
have here something new even in the case of singular cardinals of
cofinality $\omega$.

Finally, we would like to point out that the forcing construction
given in \cite{JSh} may be considered as a particular case of that
in theorem \ref{tm:d}, where the good family $\acal$ over $\lambda$
happens to be equal to the family of all proper initial segments of
$\lambda$.

\end{document}